\documentclass[12pt]{article}
\usepackage{latexsym,amsmath,amssymb}
\begin{document}
\title{On the inviscid limit for 2D incompressible flow with Navier friction condition}
\author{M.C. Lopes Filho\thanks{Research supported in part by CNPq grant \#300.962/91-6}  \\H.J. Nussenzveig Lopes\thanks{Research supported in part by CNPq grant \#300.158/93-9}\\G. Planas
\thanks{Supported by FAPESP, grant 01/14455-2}}
\date{}
\newtheorem{claim}{Claim}
\newtheorem{lemma}{Lemma}
\newtheorem{Prop}{Proposition}
\newtheorem{Theorem}{Theorem}
\newtheorem{definition}{Definition}
\newenvironment{proof}{\smallskip \noindent{\bf Proof}:}{      
               \hfill\rule{2mm}{3mm}\hspace{1in}\smallskip}
\newcommand{\disp}{\displaystyle}
\newcommand{\R}{\Bbb{R}}
\newcommand{\C}{\Bbb{C}}
\newcommand{\loc}{{\scriptsize \mbox{loc}}}
\newcommand{\vare}{\varepsilon}
\maketitle

\begin{abstract}
In \cite{Clopeau}, T. Clopeau, A. Mikeli\'c, and R. Robert studied the inviscid limit of the 2D incompressible
Navier-Stokes equations in a bounded domain subject to Navier friction-type boundary conditions. They proved that
the inviscid limit satisfies the incompressible Euler equations and their result ultimately includes flows generated 
by bounded initial vorticities. Our purpose in this article is to adapt and, to some extent, simplify their 
argument in order to include $p$-th power integrable initial vorticities, with $p>2$.   
\end{abstract}

\section{Introduction}

	In a recent paper \cite{Clopeau}, T. Clopeau, A. Mikeli\'c and R. Robert studied the inviscid limit of solutions of the 2D incompressible Navier-Stokes equations in a bounded domain with Navier friction type boundary conditions. They proved that the inviscid limit is a weak solution of the Euler equations, and their results include flows generated by bounded initial vorticities. The purpose of the present work is to extend their argument in order to include flows with initial vorticities in $L^p$, $p>2$. Technically this work involves much the same tools that were used in \cite{Clopeau} and relies in an essential manner on Clopeau {\it et alli}'s smooth data result.  

	The main motivation for studying the vanishing viscosity limit for incompressible 2D flow is the problem of boundary layers. This motivation, together with the issue of the physical meaning of the Navier friction condition was well explored in the introduction to \cite{Clopeau}.  We will not repeat that discussion here, referring the reader to that article and the references there contained for this part of our introduction. However, 2D boundary layers have been a very active field of inquiry recently, so in addition to \cite{Clopeau} we also refer the reader to \cite{dg03,tw02} for more recent developments. Beyond these issues, there is additional background which is specifically related to irregular flows which we must address here. 

	Existence of weak solutions to the incompressible 2D Euler equations has been established for rather singular initial data, more precisely, initial velocities in $L^2_{\loc}$ such that the corresponding vorticity lies in ${\cal BM}^+_c + L^1_c$, i.e. nonnegative bounded measures with compact support plus an arbitrary compactly supported integrable function. This result is due primarily to J.-M. Delort \cite{delort90}, we refer the reader also to \cite{VW93}. In both these papers the weak solutions are obtained by compactness arguments in which the initial data is mollified and the equations are subsequently exactly solved with this smooth data. Uniqueness has only been established if the initial vorticity is bounded or nearly so \cite{vishik99,yudovich63,yudovich95}, so that the issue of selection principles for singular flows is wide open. It makes sense in this case to investigate whether other approximation schemes also yield weak solutions. For example, this has been established for certain numerical schemes, see \cite{LX95,schochet96}. It is natural, from a physical point of view, to investigate the vanishing viscosity limit as well. It is possible to adapt Delort's arguments to study the inviscid limit in the absence of boundaries and this has been done for full plane flow, see \cite{majda93}. The problem of studying the existence of viscosity solutions in domains with boundary runs into the classical problem of boundary layers if one supplements the viscous approximations with the no-slip boundary condition. The work of Clopeau {\it et alli} shows that the boundary layer arising from the inviscid limit with Navier friction condition can be treated, while retaining some physical meaning. In fact, the Navier friction condition still allows for vorticity production at the boundary, but in a controlled fashion. It is therefore natural to investigate the existence of viscosity solutions by considering viscous approximations satisfying the Navier friction condition, searching for critical regularity on the initial data that guarantees the existence of such solutions. This is the main point behind the present work. 

	The remainder of this article is divided in five sections: the next section contains the basic notation and set up
of the problem; the third section investigates approximation of initial data that satisfy the Navier friction condition; the fourth section contains the {\it a priori} estimate on $L^p$-norm of vorticity which is the heart of this work; the fifth section contains a well-posedness result for the viscous approximations with $L^p$ initial vorticity; the last section contains the passage to the inviscid limit and conclusions.

\section{Preliminaries}

Let $\Omega \subseteq \R^2$ denote a bounded simply connected domain with smooth boundary. Our point of departure are the 
incompressible Navier-Stokes equations in $\Omega$. We are interested in the initial-boundary value 
problem where the velocity satisfies the Navier friction condition with friction coefficient 
$\alpha = \alpha(x) \in C^2(\partial \Omega)$, $\alpha \geq 0$. 
More precisely, the initial-boundary value problem is given by   
 
\begin{equation} \label{NSeq}
\left\{ \begin{array}{ll}  
u_t + u \cdot \nabla u = -\nabla p + \nu \Delta u & \mbox{ in } \Omega \times (0,T), \\
\mbox{div } u = 0 & \mbox{ in }\Omega\times [0,T), \\
u \cdot {\mathbf n} = 0 & \mbox{ on } \partial \Omega \times [0,T),\\
2 (Du)_S{\mathbf n} \cdot {\mathbf \tau} + \alpha u \cdot {\mathbf \tau} = 0 & \mbox{ on } \partial \Omega \times (0,T), \\
u(x,0) = u_0(x) & \mbox{ in } \Omega,  
\end{array} \right. 
\end{equation}
where $ \nu > 0 $ is the viscosity, ${\mathbf n}$ and ${\mathbf \tau}$ are the unit outwards normal and counterclockwise tangent vectors to $\partial \Omega$ respectively, $u$ is the fluid velocity, $p$ is the scalar pressure and  $(Du)_S$ is the symmetric part of the Jacobian matrix of $u$, i.e. $(Du)_S = \frac{1}{2}(Du + (Du)^t)$.  
   
The well-posedness of this initial-boundary value problem was established by 
Clopeau, Mikeli\'{c} and Robert in \cite{Clopeau}. More precisely, given a divergence-free initial velocity 
field $u_0 \in H^2(\Omega)$, tangent to the boundary, and satisfying the Navier friction condition 
$2 (Du)_S{\mathbf n} \cdot {\mathbf \tau} + \alpha u \cdot {\mathbf \tau} = 0$ on $\partial \Omega$ in the trace sense,
they showed that there exists a unique weak solution
$u^{\nu} \in L^2((0,T) ; H^1(\Omega)) \cap L^{\infty}((0,T);L^2(\Omega))$ satisfying

\begin{equation} \label{weakform} 
\begin{array}{l} \disp{
\frac{d}{dt}\int_{\Omega}  \varphi u^{\nu} + \int_{\Omega} \varphi\cdot(u^{\nu}\cdot \nabla) u^{\nu} dx \, +} 
\\ \\ \disp{
2 \nu \int_{\Omega}  (D\varphi)_S:(Du^{\nu})_S dx   
+ \nu \int_{\partial \Omega} \alpha (\varphi \cdot \tau)( u^{\nu}\cdot \tau) dS = 0,}  
\end{array} \end{equation}
for every divergence-free test vector field $\varphi \in H^1(\Omega)$, tangent to $\partial \Omega$. Here the
matrix product $A:B$ means $\sum_{i,j} A_{ij} B_{ij}$ and it is called the trace product.  
  
We note that the initial condition is not included in this weak formulation. In fact, Clopeau {\it et alli} also showed
that $u^{\nu}_t \in L^2((0,T);H^1(\Omega))$ from which it follows by integration that $u^{\nu} \in C([0,T];H^1(\Omega))$. Therefore the initial condition $u^{\nu}(\cdot,0)=u_0$ can be meaningfully imposed. Furthermore, if one assumes that the initial vorticity $\mbox{curl }u_0$ is bounded then $u^{\nu} \in C([0,T);H^2(\Omega))$. 

The Navier friction condition can be formulated in terms of vorticity. In order to do so,
a calculus identity was established in \cite{Clopeau} which we reproduce in the Lemma below.

\begin{lemma} \label{lemmaaux} Let $ v \in H^2(\Omega)$ be a vector field which is tangent to $\partial \Omega$. 
Then:
\[ (Dv)_S{\mathbf n} \cdot \mathbf{\tau} - \frac{\omega}{2}  + \kappa (v \cdot {\mathbf \tau}) = 0
\quad \mbox{ on } \partial \Omega, \]
where $ \omega = \mbox{ curl }v$ and $ \kappa $ is the curvature of $
\partial \Omega .$
\end{lemma}

One of the main difficulties in addressing the classical vanishing viscosity limit in domains with 
boundary resides in writing useful boundary conditions for the vorticity formulation of the 2D Navier-Stokes equations. 
It is through the use of the vorticity formulation that one finds higher order estimates for velocity that are 
independent of viscosity. The inviscid limit for 2D Navier-Stokes with friction condition is more tractable than the classical problem precisely because the friction boundary condition translates into a useful boundary condition for
vorticity. We introduce $\omega_0 = \mbox{ curl }u_0$ the initial vorticity and $\omega^{\nu} = \mbox{ curl } u^{\nu}$, 
the time-dependent vorticity associated to the weak solution $u^{\nu}$ of (\ref{NSeq}) with initial data $u_0$. 
For each fixed time, the velocity $u^{\nu}$ can be recovered from vorticity by means of the Biot-Savart law. We make this
explicit by writing
\[u^{\nu} = K_{\Omega}(\omega^{\nu}), \]
where $K_{\Omega}$ is an integral operator of order $-1$, with kernel given by $\nabla^{\perp}G_{\Omega}$, 
where $G_{\Omega}$ is the Green's function for the Dirichlet Laplacian in $\Omega$. Using Lemma \ref{lemmaaux} above, 
it is a standard calculation to show that $\omega^{\nu}$, $u^{\nu}$ satisfies, in a weak sense, the following parabolic initial-boundary value problem, which is the vorticity formulation of (\ref{NSeq}):

\begin{equation} \label{NSvorteq}
\left\{ \begin{array}{ll}
\omega^{\nu}_t + u^{\nu} \cdot \nabla \omega^{\nu} = \nu \Delta \omega^{\nu} & \mbox{ in } \Omega \times (0,T), \\
u^{\nu} = K_{\Omega}[\omega^{\nu}] & \mbox{ in } \Omega \times [0,T), \\
\omega^{\nu} = (2\kappa - \alpha)u^{\nu} \cdot {\mathbf \tau} & \mbox{ on } \partial \Omega \times [0,T) \\
\omega^{\nu}(\cdot,0) = \omega_0 & \mbox{ on } \Omega \times \{t=0\}.  
\end{array} \right.
\end{equation} 

\section{Approximating nonsmooth initial data}

	The problem we wish to address in this article is the inviscid limit for (\ref{NSeq}) with initial velocity 
$u_0 = K_{\Omega}[\omega_0]$, and $\omega_0 \in L^p(\Omega)$ for some $p>2$. We must first discuss this initial-boundary value problem for fixed viscosity, as this initial condition does not satisfy the conditions for the well-posedness mentioned in the previous section. It can be easily seen that this initial velocity $u_0$ is divergence free, tangent to the boundary and it belongs to $W^{1,p}(\Omega)$ (by elliptic regularity, see \cite{nu97}). This means that there is not enough regularity to impose the Navier friction condition on the initial data, so that this initial-boundary value problem is subject to an initial layer.

\begin{definition} \label{compat} We will call a function $\omega \in H^1(\Omega)\cap L^{\infty}(\Omega)$ compatible if the associated velocity $u = K_{\Omega}[\omega] \in H^2(\Omega)$ satisfies the Navier condition $\omega= (2\kappa - \alpha)u \cdot \tau$ on the boundary in the trace sense.  
\end{definition}

The first issue we need to address is how to approximate an arbitrary function in $L^p(\Omega)$ by compatible functions. This issue was addressed by Clopeau {\it et alli} for $\omega \in L^{\infty}(\Omega)$, using a fixed point argument. We will state and prove an extension of their result that applies to functions $\omega \in L^p(\Omega)$, $p>1$. The proof is a reasonably straightforward adaptation of their argument, which we include for the sake of completeness. 

\begin{lemma} \label{approx}
Let $\omega \in L^p(\Omega)$, for some $p>1$. Then there exists a sequence $\{\omega_n\}$ of compatible functions which converges to $\omega$ strongly in $L^p$.  
\end{lemma}

\begin{proof}
Recall the notation introduced in the proof of Lemma 4.2 of \cite{Clopeau}. For $x \in \overline{\Omega}$, let $d=d(x)$ be the distance of $x$ to $\partial \Omega$ and let $U_n \equiv \{x \in \overline{\Omega}: d(x) < 1/n\}$. Let $r = r(x)$ denote the orthogonal projection of $U_n$ onto $\partial \Omega$, defined for $n$ sufficiently large. Let $\zeta_n$ be a smooth cutoff for a neighborhood of $\Omega \setminus U_n$, so that $\zeta_n \equiv 0$ in $U_{n+1}$ and $\zeta_n \equiv 1$ outside $U_{n}$. Let $\eta_n$ be a standard Friedrichs mollifier. As in Lemma 4.2 we extend $\omega$ to vanish outside of $\Omega$. First, assume that $p<2$ and let $\widehat{p} = p/(2-p)$. For any $G \in L^{\widehat{p}}(\partial \Omega)$ set:
\begin{equation} \label{beta}
\beta \equiv \zeta_n (x) \eta_n \ast \omega (x) + (1-\zeta_n(x)) e^{-nd(x)} G(r(x)),
\end{equation}
and
\[v = K_{\Omega}[\beta].\]
Introduce 
\[\Psi(G) = (2\kappa - \alpha)v\cdot\tau. \]
We note that $\Psi$ maps $L^{\widehat{p}}(\partial \Omega)$ into itself. To see this, we begin by observing that $\beta \in L^p(\Omega)$. This follows since $G \in L^{\widehat{p}}(\partial \Omega)$, which implies, by a simple change of variables, that $G(r(\cdot)) \in L^{\widehat{p}}(U_n)$. As $G(r(\cdot))$ appears multiplied by a function which vanishes outside $U_n$ we may assume that $G(r(\cdot))$ vanishes outside $U_n$. Since $\widehat{p} > p$, because $p>1$, we obtain $\beta \in L^p(\Omega)$. Therefore $v \in W^{1,p}(\Omega)$, so that $v \cdot \tau \in W^{1-1/p,p}(\partial \Omega)$. We conclude using the Sobolev imbedding $W^{1-1/p,p}(\partial \Omega) \subset L^{\widehat{p}}(\partial \Omega)$.

Next we show that $\Psi$ is a contraction mapping if $n$ is sufficiently large. Let $G_1$, $G_2 \in  L^{\widehat{p}}(\partial \Omega)$. Then:
\[ \| \Psi(G_1) - \Psi(G_2)\|_{L^{\widehat{p}}(\partial\Omega)} \leq \|2\kappa - \alpha\|_{L^{\infty}}\|v_1 - v_2\|_{L^{\widehat{p}}(\partial\Omega)}
\leq C_p \|\beta_1 - \beta_2\|_{L^p(\Omega)} \]
\[ \leq C_p \|G_1(r(\cdot)) - G_2(r(\cdot))\|_{L^p(U_n)} \leq 
C_p \frac{1}{n^{1/p}}  \|G_1 - G_2\|_{L^{\widehat{p}}(\partial\Omega)}.\]
Therefore, for $n$ sufficiently large, $\Psi$ has a unique fixed point, which we denote by $G^n \in L^{\widehat{p}}(\partial \Omega)$. We denote the corresponding $\beta$ by $\omega_n$. We need to verify that $\omega_n$ is compatible. First, a standard bootstrap argument on identity (\ref{beta}), involving Sobolev imbeddings and elliptic regularity gains $1-1/\widehat{p}$ derivatives on $\omega_n$ at each step. Therefore, after a finite number of steps we reach $\omega_n \in H^1 \cap L^{\infty}$. Second, the fact that $G^n$   
is a fixed point for $\Psi$ implies that $\omega_n$ satisfies the Navier friction condition. 

Finally, we argue that $\omega_n$ converges strongly to $\omega$ in $L^p$. Since the first term on the r.h.s of (\ref{beta}) clearly converges strongly to $\omega$ in $L^p$, all we need to show is that the remaining term converges to zero in $L^p$.
First note that
\[\|(1-\zeta_n) e^{-nd(x)} G^n(r(x))\|_{L^p(\Omega)} \leq \|G^n(r(x))\|_{L^p(U_n)} \leq \frac{C}{n^{1/p}}\|G^n\|_{L^p(\partial \Omega)} \]
\[\leq {\mathit o}(1) \|G^n\|_{L^{\widehat{p}}(\partial \Omega)}.\]
Now we estimate $\|G^n\|_{L^{\widehat{p}}(\partial \Omega)}$:
\[\|G^n\|_{L^{\widehat{p}}(\partial \Omega)} \leq \|2\kappa-\alpha\|_{L^{\infty}}\|K_{\Omega}[\omega_n]\|_{L^{\widehat{p}}(\partial\Omega)} \]
\[\leq C(p,\alpha,\kappa) \|\omega_n\|_{L^p(\Omega)} \leq C (\|\omega\|_{L^p(\Omega)}
+ \|G^n(r(x))\|_{L^p(U_n)})\] 
\[ \leq C_p \|\omega\|_{L^p(\Omega)} + \frac{1}{2}\|G^n\|_{L^{\widehat{p}}(\partial\Omega)},\]
for $n$ sufficiently large, which implies the required bound.

For $p=2$ one repeats the argument above with an arbitrary $\widehat{p}$, and for $p>2$ one just
takes $\widehat{p}=\infty$. 

\end{proof}

{\bf Remark 1:} The result presented is actually more general than what we require. It applies to initial vorticities in $L^p$, $p>1$, when we are only going to use it for $p>2$. It is worth remarking that it is only for the cases $1<p\leq 2$ that we needed to use a fixed point argument in $L^{\widehat{p}}$. We could have written an argument that works for the case $p>2$ using the fixed point argument in $L^{\infty}$, like Clopeau {\it et alli} did in \cite{Clopeau}, and the proof would really be a very minor adaptation of the proof in \cite{Clopeau}, not deserving repetition even for the sake of completeness. One of the points of the present work is to clarify the criticality of this problem. This is the main reason to present the approximation result in this generality. The way it is formulated implies that this approximation issue is not part of the $p>2$ limitation. 

{\bf Remark 2:} There is no asymptotic description of the structure of the boundary layer for the present problem that would be the adaptation of the L. Prandtl description for the classical boundary layer. In the absence of such an account, the proof above gives at least a clue as to the nature of this boundary layer, embodied in the structure of the correction term. One key issue in the classical boundary layer, is that such a correction term would have, at best, an uniform $L^1$ estimate, leading to a boundary vortex sheet perturbation in the limit. This is apparent in the explicitly computable flow
generated by an impulsively started plate, known as the Rayleigh Problem, see for example \cite{panton}. This vortex sheet at the boundary is present in the inviscid limit even for smooth initial vorticities. Now, vortex sheet type regularity is critical for passing to the weak limit in approximations of the incompressible 2D Euler equations, see \cite{delort90}. In some sense, it is this fact that is the heart of the difficulty in the classical boundary layer problem. The correction term in the proof above suggests that the boundary layer associated to the Navier friction condition would correspond to uniformly bounded vorticity near the boundary for $p>2$, and $L^{\widehat p}$ vorticity near the boundary for $1<p\leq 2$, so that one would expect criticality only at $p=1$.

{\bf Remark 3:} The argument presented breaks down when $p=1$, mainly because elliptic regularity breaks down, so that one cannot guarantee that $\Psi$ maps $L^1$ to itself.

\section{{\it A priori} estimate on vorticity}

The purpose of this section is to derive an {\it a priori} estimate for vorticity on solutions of (\ref{NSeq}). We begin with a compatible initial vorticity $\omega_0$, as defined in the previous section. We use $u_0 = K_{\Omega}[\omega_0]$ as initial data. The well-posedness of the initial-boundary value problem (\ref{NSeq}) for such initial data was established in \cite{Clopeau}, as previously mentioned. Let $u = u(x,t)$ be the unique weak solution of (\ref{NSeq}) with data $u_0$.   
The vector field $u$ belongs to $C([0,\infty);H^2(\Omega))$ and it satisfies the weak formulation (\ref{weakform}) of the Navier-Stokes equation with Navier friction condition.  The vorticity $\omega = \mbox{curl }u $ satisfies the parabolic
equation (\ref{NSvorteq}) in a weak sense.

\begin{lemma} \label{apriori}  Fix $p>2$. There exists a constant $C>0$, depending only on $p$, $\Omega$ and the friction coefficient $\alpha$ such that the vorticity satisfies:
\[\|\omega(\cdot,t)\|_{L^p} \leq   C (\|\omega_0\|_{L^p} + \|u_0\|_{L^2}). \] 
\end{lemma}
 
\begin{proof}
The proof involves applying a maximum principle to two auxiliary problems. First observe that $u\cdot \tau \in L^{\infty}(\partial \Omega \times (0,T))$ since $u \in C([0,T); H^2(\Omega))$. Set 
\[\Lambda = \|(2\kappa - \alpha)u\cdot\tau\|_{L^{\infty}(\partial \Omega \times (0,T))}.\] 
Consider the initial-boundary value problem for the Fokker-Planck equation:
 
\begin{equation} \label{FP1}
\left \{ \begin{array}{ll}
 \widetilde{\omega}_t - \nu \Delta \widetilde{\omega} + u \cdot
\nabla \widetilde{\omega}=0 & \mbox{ in } \Omega \times (0,T),\\
\widetilde{\omega}(\cdot,0) = |\omega_0| & \mbox { in } \Omega,\\
\widetilde{\omega} = \Lambda & \mbox{ on }
\partial \Omega \times (0,T).
\end{array} \right.
\end{equation}

This problem has a unique weak solution $\widetilde{\omega} \in L^2((0,T);H^1(\Omega))$, by
Theorem 6.1 and 6.2 in \cite{lieberman}. Then, $ \omega_1 = \omega - \widetilde{\omega}$ is a weak solution 
for the following initial-boundary value problem:

\begin{equation} \label{FP2} \left \{ \begin{array}{ll}
(\omega_1)_t - \nu \Delta \omega_1 + u \cdot
\nabla \omega_1 = 0 & \mbox{ in } \Omega \times (0,T),\\
\omega_1 (\cdot,0) = \omega_0 - |\omega_0| & \mbox { in } \Omega,\\
\omega_1 = ( 2 \kappa - \alpha ) u \cdot \tau - \Lambda & \mbox{ on }
\partial \Omega \times (0,T).
\end{array} \right.
\end{equation}
The coefficients of the Fokker-Planck operator $\partial_t -\nu \Delta + u \cdot \nabla$ are such that the maximum principle for weak solutions, given in Corollary 6.26 of \cite{lieberman}, is valid. Therefore, as $\omega_1 \leq 0$ on the parabolic boundary $\partial \Omega \times (0,T)\cup \Omega \times \{t=0\}$, it follows that $\omega_1 \leq 0$ a.e. in $\Omega \times [0,T)$. Analogously, we prove that $\omega_2 = -\omega - \widetilde{\omega} $ is non-positive. We thus obtain

\begin{equation} \label{estimativa}
| \omega | \leq \widetilde{\omega} \mbox{ a.e. in } \Omega \times [0,T).
\end{equation}

Moreover, as $\omega_0$ is compatible it is bounded. Hence Corollary 6.26 of \cite{lieberman} may also be applied to equation (\ref{FP1}) yielding that $\widetilde{\omega} \in L^{\infty}((0,T)\times \Omega)$.

Next we obtain an estimate for $ \widetilde{\omega}. $ Let $
\widehat{\omega} = \widetilde{\omega} - \Lambda $. This is a solution of the following
problem:

\begin{equation} \label{widehateq}
 \left \{ \begin{array}{ll}
 \widehat{\omega}_t - \nu \Delta \widehat{\omega} + u \cdot
\nabla \widehat{\omega}=0 & \mbox{ in } \Omega \times (0,T),\\
\widehat{\omega}(\cdot,0) = |\omega_0| - \Lambda & \mbox { in } \Omega,\\
\widehat{\omega} = 0 & \mbox{ on }
\partial \Omega \times (0,T).
\end{array} \right.
\end{equation}

We formally multiply (\ref{widehateq}) by $ \widehat{\omega} |\widehat{\omega}|^{p-2}$, where $ p > 2$, we integrate by parts and use the incompressibility of the flow $u$ to obtain:
\[ \frac{1}{p} \frac{d}{dt} \int_{\Omega} | \widehat{\omega} |^p +
(p-1) \nu \int_\Omega  ||\nabla \widehat{\omega}
||\widehat{\omega}|^{(p-2)/2} | ^2 dx = 0. \] Then,
\[
\| \widehat{\omega}(\cdot,t)\|_{L^p(\Omega)}  \leq  \|
\widehat{\omega}(\cdot,0)\|_{L^p(\Omega)}  \leq  \|\omega_0 \|_{L^p(\Omega)}
+ \Lambda |\Omega |^{1/p}. \] Therefore,
\begin{equation} \label{widetildest}
\| \widetilde{\omega} \|_{L^p(\Omega)} \leq \|
\widehat{\omega}\|_{L^p(\Omega)} + \Lambda |\Omega |^{1/p} \leq \|\omega_0
\|_{L^p(\Omega)} + 2 \Lambda |\Omega |^{1/p}.
\end{equation}

This formal calculation can be made rigorous by using the weak formulation of (\ref{widehateq}) given in \cite{lieberman}. One needs to approximate the function $\widehat{\omega} |\widehat{\omega}|^{p-2}$ by suitable smooth test functions in such a way as to pass to the limit in each term of the weak formulation. The relevant observations in this justification are that $\widehat{\omega}_t \in L^2((0,T);H^{-1}(\Omega))$ and $\widehat{\omega} \in L^2((0,T);H^1_0(\Omega)) \cap L^{\infty}((0,T)\times \Omega )$, so that the following identity holds in 
$L^2((0,T);H^{-1}(\Omega))$:
\[\partial_t\left( \frac{1}{p}|\widehat{\omega}|^p\right) = |\widehat{\omega}|^{p-2}\widehat{\omega}\widehat{\omega}_t . \]

Given (\ref{widetildest}) we now turn to the estimate of $\Lambda$. Using Sobolev imbedding and interpolating between $W^{1,p}$ and $L^2$, we find:
\[\| u (\cdot,t) \cdot \tau \|_{L^{\infty}(\partial \Omega)} \leq C \|u (\cdot,t) \|_{C(\bar{\Omega})} 
\leq C \|u (\cdot,t)\|_{L^2(\Omega)}^\theta \| u (\cdot,t) \|_{W^{1,p}(\Omega)}^{1- \theta} \]
\[\leq C \|u (\cdot,t) \|_{L^2(\Omega)}^\theta \| \omega (\cdot,t)\|_{L^p(\Omega)}^{1- \theta},\] 
where $\theta = (p-2)/(2p-2)$.
 
Let $\varepsilon$ be an arbitrary positive number. We now use Young's inequality together with the fact that $\kappa$ and $\alpha$ are bounded to conclude that:
\begin{equation} \label{estimativa2}
\Lambda \leq C_{\varepsilon}  \|u  \|_{L^{\infty}((0,T);L^2(\Omega)} +  \varepsilon \| \omega \|_{L^{\infty}((0,T);L^p(\Omega))} \end{equation}
for some $C_{\varepsilon}>0$. Taking $\varepsilon$ small enough, from
(\ref{estimativa})-(\ref{estimativa2}) we obtain:
\begin{equation} \label{estimatevort}
\| \omega \|_{L^\infty(0,T;L^p(\Omega))}  \leq C(
\|\omega_0 \|_{L^p(\Omega)}  + \|u \|_{L^\infty(0,T;L^2(\Omega))})  \end{equation} 
for any $p > 2$, where $ C = C(p,\Omega,\|\kappa\|_{L^{\infty}(\partial
\Omega)},\|\alpha\|_{L^\infty(\partial \Omega)} ).$ Finally, a standard energy estimate, such as the one carried out in \cite{Clopeau} (see estimate (2.16)), yields 
$\|u \|_{L^\infty(0,T;L^2(\Omega))} \leq \|u_0\|_{L^2(\Omega)}$, thereby concluding the proof. 

\end{proof}

{\bf Remark 1:} This Lemma is the heart of this article. Note that the restriction $p>2$ comes into the proof above 
because of the need to produce an uniform bound on the velocity at the boundary. It would be interesting to know if this is a physically meaningful restriction. This would mean that the problem of controlling the generation of vorticity by the interaction of incompressible flow with a "Navier condition" boundary is critical at $p=2$. 
However, this criticality at $p=2$ seems unlikely. The limitation on the integrability of vorticity in the proof above appears to reflect a limitation on the maximum principle technique employed rather than an essential feature of this problem. In contrast, the exponent $p=1$ found to be critical in the proof of Lemma \ref{approx} seems much more essential and is already known to be critical in terms of passage to weak limits on the nonlinearity of the incompressible 2D Euler and Navier-Stokes equations. 

{\bf Remark 2:} The natural way to extend this vorticity estimate to $p \leq 2$ would be to derive an $L^p$ energy estimate on the vorticity equation. Multiplying the vorticity equation (\ref{NSvorteq}) by $p\omega|\omega|^{p-2}$, integrating in space and performing the usual integration by parts yields: 
\[\frac{d}{dt} \int_{\Omega} |\omega|^p dx = -\nu p (p-1) \int_{\Omega} |\omega|^{p-2} |\nabla \omega|^2 dx + \nu p \int_{\partial \Omega}|\omega|^{p-2} \omega \nabla \omega \cdot \mathbf{n} dS. \]  
We note that the boundary term is the flux of $|\omega|^p$ through the boundary, over which we have no control. One special case for which this simple estimate does provide an improvement over Lemma \ref{apriori} is the case of $\alpha = 2 \kappa$, as then the troublesome boundary term vanishes. This corresponds to the so called free boundary condition $\omega=0$ on $\partial \Omega$. It is a well-known fact that one can handle the inviscid asymptotics in this case, as one is imposing that the boundary does not generate vorticity and thus there are no boundary layers. For details, see \cite{jllions} and the special case of time-independent domain in \cite{hh00}.

\section{Well-posedness for the viscous problem}

In this section we observe that the initial-boundary value problem for the Navier-Stokes equations with friction-type boundary condition is well-posed even if the initial vorticity is not compatible. This was already done in \cite{Clopeau} for bounded initial vorticity. 

To begin with, we require a notion of weak solution that is weaker than (\ref{weakform}), which we will
introduce in the result below.

\begin{lemma} \label{weakweakform}
Let $u^{\nu} \in L^2((0,T) ; H^1(\Omega)) \cap L^{\infty}((0.T);L^2(\Omega))$ be a weak solution of (\ref{NSeq}) in the sense given by the identity (\ref{weakform}). Then for any test vector field $\varphi \in C^{\infty}_c([0,T)\times \overline{\Omega})$,
divergence free and tangent to $\partial \Omega$ we have:
\[\int_0^T \int_{\Omega} u^{\nu} \varphi_t + u^{\nu}(u^{\nu} \cdot \nabla )\varphi dxdt + \int_{\Omega} u_0 \varphi(\cdot,0)dx \]
\[ = 2\nu \int_0^T\int_{\Omega} (D\varphi)_S : (Du)_S dxdt + \nu \int_0^T \int_{\partial\Omega} \alpha (\varphi \cdot \tau)(u^{\nu} \cdot \tau)dS dt. \]  
\end{lemma}

\begin{proof}
Let $\varphi$ be a test vector field. For each $s \in [0,T)$, define 
\[g(t,s) \equiv \int_{\Omega} u^{\nu}(x,t) \varphi(x,s) dx.\]
Then, by (\ref{weakform}) we have:
\[\frac{\partial g}{\partial t} = -\int_{\Omega} \varphi\cdot(u^{\nu}\cdot \nabla) u^{\nu} dx \]
\[ - 2 \nu \int_{\Omega}  (D\varphi)_S:(Du^{\nu})_S dx   
- \nu \int_{\partial \Omega} \alpha (\varphi \cdot \tau)( u^{\nu}\cdot \tau) dS \]
\[= \int_{\Omega} u^{\nu}(u^{\nu}\cdot \nabla)\varphi dx - 2 \nu \int_{\Omega}  (D\varphi)_S:(Du^{\nu})_S dx   
- \nu \int_{\partial \Omega} \alpha (\varphi \cdot \tau)( u^{\nu}\cdot \tau) dS, \]
using integration by parts and the fact that $u^{\nu}$ is divergence free. On the other hand, we also have:
\[\frac{\partial g}{\partial s} = \int_{\Omega} u^{\nu}(x,t) \varphi_s(x,s) dx. \]
Therefore, it follows that:
\[\frac{d}{dt} (g(t,t)) = \int_{\Omega} u^{\nu}(x,t) \varphi_t(x,t) dx \]
\[ + \int_{\Omega} u^{\nu}(u^{\nu}\cdot \nabla)\varphi dx - 2 \nu \int_{\Omega}  (D\varphi)_S:(Du^{\nu})_S dx   
- \nu \int_{\partial \Omega} \alpha (\varphi \cdot \tau)( u^{\nu}\cdot \tau) dS.\]
Integrating this last identity in time and identifying the initial data yields the desired result.

\end{proof}

We now state and prove the main result in this section.

\begin{Prop} \label{viscwp} Let $\omega_0 \in L^p(\Omega)$ for some $p>2$ and $u_0 = K_{\Omega}[\omega_0]$. Fix $\nu > 0$. Then there exists a unique vector field $u^{\nu} = u^{\nu}(x,t) \in C([0,T);L^2(\Omega)) \cap L^2((0,T);H^1(\Omega))$ satisfying the weak formulation (\ref{weakform}) of the 2D incompressible Navier-Stokes system (\ref{NSeq}) with initial data $u_0$. Moreover, the associated vorticity $\omega^{\nu} = \mbox{ curl }u^{\nu}$ satisfies the estimate 
\[ \|\omega^{\nu}(\cdot,t)\|_{L^p(\Omega)} \leq C,\] 
a.e. in time, with constant $C>0$ independent of viscosity.
\end{Prop}

\begin{proof}
Let $\omega_{0,n}$ be a sequence of compatible functions approximating $\omega_0$ strongly in $L^p$, as constructed in
Lemma \ref{approx} and $u_{0,n} = K_{\Omega}[\omega_{0,n}]$. Let $u^{\nu}_n$ be the weak solution of system 
(\ref{NSeq}) given by the well-posedness result of \cite{Clopeau} and $\omega^{\nu}_n = \mbox{ curl }u^{\nu}_n$ the
corresponding vorticity. We begin by observing that Lemma \ref{apriori} gives the uniform estimate:
\begin{equation} \label{vrtest}
\|\omega^{\nu}_n\|_{L^{\infty}((0,T);L^p(\Omega))} \leq C (\|\omega_0\|_{L^p(\Omega)} + \|u_0\|_{L^2(\Omega)}),
\end{equation}
for some $C>0$. By the Poincar\'{e} and Calder\'{o}n-Zygmund inequalities it follows that 
\begin{equation} \label{velest0}
\|u^{\nu}_n\|_{L^{\infty}((0,T);W^{1,p}(\Omega))} \leq C (\|\omega_0\|_{L^p(\Omega)} + \|u_0\|_{L^2(\Omega)}).
\end{equation}

Let $\varphi \in C^{\infty}_c((0,T)\times\overline{\Omega})$ be a divergence free test vector field which is tangent to the boundary of $\Omega$. We compute the time-derivative of $u^{\nu}_n$ in the sense of distributions. This is where we require the new weak formulation given in Lemma \ref{weakweakform}. We have:
\[ \langle \varphi,\partial_t u^{\nu}_n \rangle = - \int_0^T \int_{\Omega}(\partial_t \varphi) u^{\nu}_n \]
\[= \int_0^T\int_{\Omega} u^{\nu}_n(u^{\nu}_n \cdot \nabla )\varphi -2\nu (D\varphi)_S : (Du^{\nu}_n)_S dxdt - 
\nu \int_0^T \int_{\partial \Omega} \alpha (\varphi \cdot \tau)( u^{\nu}_n \cdot \tau) dS.  \]  

Recall that $p>2$, so that (\ref{velest0}) implies that $\|u^{\nu}_n\|_{L^{\infty}((0,T)\times\Omega)} \leq C$, for some constant $C>0$ depending only on the initial data. Similarly, $\|Du^{\nu}_n\|_{L^2((0,T)\times\Omega)}$ is bounded uniformly by a positive constant depending only on initial data. We use these facts to estimate $\partial_t u^{\nu}_n$. Let $\varphi$ be a test vector field, which we first assume to be divergence free as above. We have:
\[|\langle \varphi,\partial_t u^{\nu}_n \rangle| \leq 
\left( \|u^{\nu}_n\|_{L^{\infty}((0,T)\times\Omega)}\|u^{\nu}_n\|_{L^2((0,T)\times\Omega)}
+ 2\nu \|Du^{\nu}_n\|_{L^2((0,T)\times\Omega)} \right.\]
\[\left. + \nu \|\alpha u^{\nu}_n\|_{L^2((0,T);H^1(\Omega))}\right) \|\varphi\|_{L^2((0,T);H^1(\Omega))} \leq C \|\varphi\|_{L^2((0,T);H^1(\Omega))}, \]
where we have used the continuity of the trace operator from $H^1(\Omega)$ onto $L^2(\partial \Omega)$ to estimate the boundary term. 
Now, if the test vector field $\varphi$ is not divergence free, we use standard properties of the Leray projector $P$ to obtain an estimate 
\[\|P\varphi\|_{L^2((0,T);H^1(\Omega))} \leq C  \|\varphi\|_{L^2((0,T);H^1(\Omega))}, \]
and we repeat the argument above with $P\varphi$ in place of $\varphi$. Note that 
\[\langle  \varphi,\partial_t u^{\nu}_n \rangle = 
\langle  P\varphi,\partial_t u^{\nu}_n \rangle\]
as $\partial_t u^{\nu}_n$ is divergence free and tangent to the boundary. By duality this implies the following estimate:
\begin{equation} \label{velest1}
\|\partial_t u^{\nu}_n\|_{L^2((0,T);H^{-1}(\Omega))} \leq C,  
\end{equation}
with $C>0$ depending only on the initial data. Thus $u^{\nu}_n$ is equicontinuous from $(0,T)$ to $H^{-1}(\Omega)$ and we can use the Aubin-Lions Lemma to obtain a subsequence, which we will not relabel, converging strongly in $C([0,T];L^2(\Omega))$. Without loss of generality this subsequence also converges weakly in $L^2((0,T);H^1(\Omega))$ to a limit $u^{\nu}$. It is now easy to see that we can pass to the limit in each term in the weak formulation (\ref{weakform}) of the Navier-Stokes equations, thereby concluding the proof of existence for the initial-boundary value problem (\ref{NSeq}). Furthermore, from the estimate on vorticity (\ref{vrtest}) it follows that 
\[\|\omega^{\nu}(\cdot,t)\|_{L^p(\Omega)} \leq C,\]  a.e. in time, for some constant $C>0$ depending only on the initial data. 

The uniqueness portion of this result is standard,  and may be obtained by adapting the classical argument using an energy estimate on the difference of two solutions with the same data. We conclude that there exists at most one weak solution $u \in C([0,T];L^2(\Omega)) \cap L^2((0,T);H^1(\Omega))$ of (\ref{NSeq}) in the form (\ref{weakform}). 

\end{proof}

{\bf Remark:} The proof above can be adapted for $1<p\leq 2$, assuming of course that Lemma \ref{apriori} could be proved
in that case. The main steps in this adaptation would be:
\begin{itemize} 
\item substitute the $L^{\infty}$ estimate on $u^{\nu}_n$ by an 
$L^{\infty}((0,T);L^{p^{\ast}})$ estimate, with $p^{\ast}$ either the critical Sobolev exponent if $p<2$ or an arbitrary number $1<q<\infty$ if $p=2$;
\item use the fact that $\sqrt{\nu} \|u^{\nu}_n\|_{L^2((0,T);H^1(\Omega))}$ is bounded uniformly in $n$ and $\nu$ by the $L^2$-norm of the initial velocity. This is a consequence of standard energy estimates for the Navier-Stokes equations. 
\end{itemize}

\section{Inviscid limit and Conclusions}

Let $\omega_0 \in L^p(\Omega)$ for some $p > 2$ and let $u_0 = K_{\Omega}[\omega_0]$. In this last section we show that the sequence of solutions of the Navier-Stokes equations with initial velocity $u_0$ and with Navier friction conditions possesses a converging subsequence to a solution of the Euler equations with same initial velocity as viscosity vanishes. 
The proof is very similar to the existence part of the proof of Proposition \ref{viscwp}.  

\begin{Theorem} \label{inviscidlim}
Let $u^{\nu} = u^{\nu}(x,t)$ be the solution of (\ref{NSeq}) such that $u^{\nu}(\cdot,0)=u_0$. Then there exists a sequence $\nu_k \to 0$ such that $u^{\nu_k} \to u$ strongly in $C([0,T];L^2(\Omega))$ as $k \to \infty$ and $u$ is a weak solution of the incompressible 2D Euler equations in the sense:
\[\int_0^T \int_{\Omega} u \varphi_t + u(u \cdot \nabla )\varphi dxdt + \int_{\Omega} u_0 \varphi(\cdot,0)dx =0,\]
for any test vector field $\varphi \in C^{\infty}_c([0,T)\times\overline{\Omega})$ which is divergence free and tangent to the boundary.   
\end{Theorem}

\begin{proof}
We recall from the proof of Proposition \ref{viscwp} that the following uniform estimates hold for $u^{\nu}$: 
\[\|u^{\nu}\|_{L^{\infty}((0,T);W^{1,p}(\Omega))} \leq C\]
and 
\[\|\partial_t u^{\nu}\|_{L^2((0,T);H^{-1}(\Omega))} \leq C,\]
where $C>0$ depends only on the initial velocity $u_0$ and initial vorticity $\omega_0$
and is independent of viscosity (see the proof of (\ref{velest0}) and (\ref{velest1})). 
From these estimates it is possible to extract a subsequence $u^{\nu_k}$ which converges strongly in $C([0,T];L^2(\Omega))$ and weakly in $L^2((0,T);H^1(\Omega))$. It is easy to see that these modes of convergence are sufficient to pass to the limit in each term of the weaker weak formulation, given in Lemma \ref{weakweakform}, and guarantee that the limit function $u$ satisfies the identity:
\[ \int_0^T \int_{\Omega} u \varphi_t + u (u \cdot \nabla )\varphi dxdt + \int_{\Omega} u_0 \varphi(\cdot,0)dx =0\]
for any test vector field $\varphi \in C^{\infty}_c([0,T)\times \overline{\Omega})$
which is divergence free and tangent to $\partial \Omega$. This is precisely the standard formulation of a weak solution of the Euler equations, hence we conclude the proof.

\end{proof}

{\bf Remark:} We could have used the weak formulation given by (\ref{weakform}) to pass to the inviscid limit, thereby
obtaining a weak solution to the Euler equations satisfying (\ref{weakform}) with $\nu = 0$. We chose to use the form in Theorem \ref{inviscidlim} because it is the standard weak formulation of the incompressible 2D Euler equations. 

We conclude this article with a few final observations. First, we call attention once more to the fact that the authors are not convinced of the criticality of $p=2$, so the critical $p$ remains an open problem. Second, as mentioned in Section 3, there is no asymptotic description in the fluid mechanics literature of the boundary layer associated with the Navier friction condition, something that, if available, would clarify the issues raised here. Finally, an interesting question which we have not explored is whether the viscosity weak solution of the incompressible 2D Euler equations obtained above conserves $L^p$-norm of vorticity. Conservation of the $L^p$ norm of vorticity holds both for weak solutions in the full plane are known to do, as a consequence of DiPerna-Lions theory, see \cite{lionsbook} and for strong solutions, as one can ascertain directly from the vorticity equation. In the viscous approximation, vorticity can be generated at the boundary, so that the question is whether this possibility disappears in the vanishing viscosity regime.      

\vspace{.5cm}

{\it Acknowledgements:} The authors would like to thank D. Iftimie and M. O. Souza for enlightening discussions.

\vskip\baselineskip
\noindent
{\sc
Milton C. Lopes Filho\\
Departamento de Matem\'{a}tica, IMECC-UNICAMP.\\
Caixa Postal 6065, Campinas, SP 13083-970, Brasil
\\}
{\it E-mail address:} mlopes@ime.unicamp.br

\vspace{.1in}
\noindent
{\sc 
Helena J. Nussenzveig Lopes\\
Departamento de Matem\'{a}tica, IMECC-UNICAMP.\\
Caixa Postal 6065, Campinas, SP 13083-970, Brasil
\\}
{\it E-mail address:} hlopes@ime.unicamp.br

\vspace{.1in}
\noindent
{\sc
Gabriela del Valle Planas\\
Departamento de Matem\'{a}tica, ICMC-USP.\\
Caixa Postal 668, S\~{a}o Carlos, SP 13560-970, Brasil
\\}
{\it E-mail address:} gplanas@icmc.usp.br


\vspace{.5cm}

\begin{thebibliography}{99}

\bibitem{Clopeau} Clopeau, T., Mikeli\'c, A., Robert, R., {\it On the
vanishing viscosity limit for the 2D incompressible Navier-Stokes
equations with the friction type boundary conditions}, 
Nonlinearity  {\bf 11} (1998) 1625--1636.

 
\bibitem{delort90} Delort, J.-M., {\it Existence de nappes de tourbillon en dimension deux}, J. Amer. Math. Soc.  {\bf 4}  (1991)  no. 3, 
553--586.

\bibitem{dg03} Desjardins, B. and Grenier, E., {\it Linear instability implies nonlinear instability for various types of viscous boundary layers},  Ann. Inst. H. Poincar\'{e} Anal. Non Lin\'{e}aire  {\bf 20} (2003), 87--106.

\bibitem{hh00} He, C. and Hsiao, L., {\it Two-dimensional Euler equations in a time dependent domain}, J. Diff. Eqns. {\bf 163} (2000) 
no. 2, 265--291.

\bibitem{lieberman} Lieberman, G., {\it  Second order parabolic differential equations.} World Scientific Publishing Co., Inc., River Edge, NJ, 1996.  

\bibitem{jllions} Lions, J.L., {\it Quelques m\'ethodes de resolution des probl\'emes aux limites non lin\'eaires}, Dunod, Gauthier-Villars, Paris, 1969.

\bibitem{lionsbook} Lions, P.L., {\it Mathematical Topics in Fluid Mechanics v. I, Incompressible Models}, Oxford Lecture Ser. in Mathematics and its Applications v. 3, Clarendon Press, Oxford, 1996.

\bibitem{LX95} Liu, J.-G. and Xin, Z., {\it Convergence of vortex methods for weak solutions to the $2$-D Euler equations with vortex sheet data},  Comm. Pure Appl. Math. {\bf 48}  (1995)   611--628.

\bibitem{majda93} Majda, A. J., {\it  Remarks on weak solutions for vortex sheets with a distinguished sign},  Indiana Univ. Math. J. {\bf 42}  (1993)  921--939. 

\bibitem{nu97}  Nussenzveig Lopes, H. J., {\it A refined estimate of the size of concentration sets for $2$D incompressible inviscid flow}, Indiana Univ. Math. J.  {\bf 46}  (1997) 165--182.

\bibitem{panton} Panton, R. L., {\it Incompressible Flow} John Wiley and Sons, New York, 1984.

\bibitem{schochet96} Schochet, S., {\it The point-vortex method for periodic weak solutions of the 2-D Euler equations}, Comm. Pure Appl. Math. {\bf 49}  (1996) 911--965. 

\bibitem{tw02} Temam, R. and Wang, X., {\it Boundary layers associated with incompressible Navier-Stokes equations: the noncharacteristic boundary case}, J. Differential Equations {\bf 179} (2002), 647--686.
 
\bibitem{VW93} Vecchi, I. and Wu, S., {\it On $L\sp 1$-vorticity for $2$-D incompressible flow},  Manuscripta Math.  
{\bf 78}  (1993)  403--412. 

\bibitem{vishik99} Vishik, M., {\it Incompressible flows of an ideal fluid with vorticity in borderline spaces of Besov type}, Ann. Sci. Ecole Norm. Sup. (4) {\bf 32} (1999) 769--812.

\bibitem{yudovich63} Yudovich, V. I., {\it Non-stationary flows of an ideal incompressible fluid}, (Russian)
\u Z. Vy\v cisl. Mat. i Mat. Fiz. {\bf 3} (1963) 1032--1066.

\bibitem{yudovich95} Yudovich, V. I., {\it Uniqueness theorem for the basic nonstationary problem in the dynamics of an ideal incompressible fluid},  Math. Res. Lett.  {\bf 2}  (1995)  27--38. 

\end{thebibliography}
\end{document}